\RequirePackage{ifpdf}
\ifpdf 
\documentclass[pdftex]{sigma}
\else
\documentclass{sigma}
\fi

\begin{document}
\allowdisplaybreaks

\renewcommand{\PaperNumber}{080}

\renewcommand{\thefootnote}{$\star$}

\FirstPageHeading

\ShortArticleName{A Formula for the Logarithm of the KZ Associator}

\ArticleName{A Formula for the Logarithm of the KZ Associator\footnote{This paper is a contribution 
to the Vadim Kuznetsov Memorial Issue ``Integrable Systems and Related Topics''.
The full collection is available at 
\href{http://www.emis.de/journals/SIGMA/kuznetsov.html}{http://www.emis.de/journals/SIGMA/kuznetsov.html}}}

\Author{Benjamin ENRIQUEZ~$^\dag$ and Fabio GAVARINI~$^\ddag$}
\AuthorNameForHeading{B. Enriquez and F. Gavarini}

\Address{$^\dag$~IRMA (CNRS), rue Ren\'e Descartes, F-67084 Strasbourg, France} 

\EmailD{\href{mailto:enriquez@@math.u-strasbg.fr}{enriquez@@math.u-strasbg.fr}} 

\Address{$^\ddag$~Universit\'a degli Studi di Roma ``Tor Vergata'', Dipartimento di 
Matematica,\\ 
$\phantom{^\ddag}$~Via della Ricerca Scientif\/ica 1, I-00133 Rome, Italy}

\EmailD{\href{mailto:gavarini@@mat.uniroma2.it}{gavarini@@mat.uniroma2.it}}

\ArticleDates{Received October 03, 2006, in f\/inal form November 10,
2006; Published online November 13, 2006}

\Abstract{We prove that the logarithm of a group-like element in a free algebra
coincides with its image by a certain linear map.  We use this result
and the formula of Le and Murakami for the Knizhnik--Zamolodchikov (KZ)
associator $\Phi$ to derive a formula for $\log(\Phi)$ in terms of
MZV's (multiple zeta values).}


\Classification{17B01; 81R50} 

\begin{flushright}
\it To the memory of Vadim Kuznetsov.
\end{flushright}

\section{Logarithms of group-like elements}

Let $F_n$ be the free associative algebra generated by free variables
$x_1,\dots,x_n$, let ${\mathfrak f}_n\subset F_n$ be the free Lie algebra with the
same generators, and let $\widehat{\mathfrak f}_n$, $\widehat F_n$ be their degree
completions (where $x_1,\dots,x_n$ have degree $1$). A group-like
element of $\widehat F_n$ is an element of the form $X = 1 +$ (terms of
degree $>0)$, such that $\Delta(X) = X\otimes X$, where $\Delta$ is
the completion of the coproduct $F_n \to F_n^{\otimes 2}$, for which
$x_1,\dots,x_n$ are primitive. It is well-known that the exponential
def\/ines a bijection $\exp$:~$\widehat{\mathfrak f}_n \to \{$group-like elements of
$\widehat F_n\}$ (also denoted $x\mapsto e^x$).  We denote by log the
inverse bijection.

We denote by ${\rm CBH}_n(x_1,\dots ,x_n)$ the multilinear part (in
$x_1,\dots ,x_n$) of $\log (e^{x_1}\cdots e^{x_n})$. Def\/ine a linear map
\[
{\rm cbh}_n : \ F_n \to {\mathfrak f}_n 
\] 
by ${\rm cbh}_n(1)=0$ and ${\rm cbh}_n(x_{i_1}\cdots x_{i_k}):=
{\rm CBH}_k(x_{i_1},\dots ,x_{i_k})$. This map extends to 
a linear map $\widehat{{\rm cbh}}_n$: $\widehat F_n \to \widehat {\mathfrak f}_n$. 

\begin{proposition} \label{main}
If $X\in \widehat F_n$ is group-like, then $\log (X) = \widehat{{\rm cbh}}_n(X)$.
\end{proposition}

\begin{proof} It is known that $F_n = U({\mathfrak f}_n)$, so that the
symmetrization is an isomorphism \linebreak ${\rm sym}$:~$S({\mathfrak f}_n) \to F_n$.  Denote
by $p_n$: $F_n \to {\mathfrak f}_n$ the composition of ${\rm sym}^{-1}$ with the
projection $S({\mathfrak f}_n) = \oplus_{k\geq 0} S^k({\mathfrak f}_n)$ onto $S^1({\mathfrak f}_n) =
{\mathfrak f}_n$. We f\/irst prove: 

\begin{lemma}
$p_n = {\rm cbh}_n$. 
\end{lemma}

\begin{proof} If ${\mathfrak g}$ is a Lie algebra, let $p_{\mathfrak g}$: $U({\mathfrak g})\to{\mathfrak g}$
be the composition of the inverse of the symmetrization 
$S({\mathfrak g})\to U({\mathfrak g})$ with the projection onto the f\/irst component
of $S({\mathfrak g})$. If $\phi$: ${\mathfrak g}\to{\mathfrak h}$ is a~Lie algebra morphism, then we 
have a commutative diagram 
\[
\begin{matrix}
U({\mathfrak g}) & \stackrel{p_{\mathfrak g}}{\to}& {\mathfrak g} \\
\scriptstyle{U(\phi)}\downarrow & & \downarrow\scriptstyle{\phi} \\
U({\mathfrak h}) & \stackrel{p_{\mathfrak h}}{\to}& {\mathfrak h} 
\end{matrix}
\]
It follows from Lemma 3.10 of \cite{R} that if 
$k\geq 0$ and $F_k$ is the free algebra with generators 
$y_1,\dots ,y_k$, then $p_k(y_1,\dots ,y_k) = {\rm CBH}_k(y_1,\dots ,y_k)$.
If now ${\boldsymbol i} = (i_1,\dots ,i_k)$ is a sequence of elements of 
$\{1,\dots ,n\}$, 
we have a unique morphism $\phi_{{\boldsymbol i}}$: ${\mathfrak f}_k\to {\mathfrak f}_n$, 
such that $y_1\mapsto x_{i_1}$, \dots, $y_k\mapsto x_{i_k}$.

Then 
\begin{gather*}
 p_n(x_{i_1}\cdots x_{i_k}) = p_{{\mathfrak f}_n} \circ U(\phi_{{\boldsymbol i}})(y_1\cdots y_k)
= \phi_{{\boldsymbol i}} \circ p_{{\mathfrak f}_k}(y_{1}\cdots y_{k}) 
\\ 
\phantom{p_n(x_{i_1}\cdots x_{i_k})}{} = \phi_{{\boldsymbol i}}
({\rm CBH}_k(y_1,\dots ,y_k)) = {\rm CBH}_k(x_{i_1},\dots ,x_{i_k}) = 
{\rm cbh}_n(x_{i_1}\cdots x_{i_k}), 
\end{gather*} 
which proves the lemma. 
\end{proof}

\noindent
{\bf End of proof of Proposition \ref{main}.}\
We denote by $\widehat p_n : \widehat F_n \to \widehat{\mathfrak f}_n$ the map similarly derived
from the isomorphism $\widehat F_n \simeq \widehat\oplus_{k\geq 0} S^k(\widehat{\mathfrak f}_n)$
(where $\widehat\oplus$ is the direct product). Then $p_n={\rm cbh}_n$ implies
$\widehat p_n = \widehat{{\rm cbh}}_n$.

If now $X \in \widehat F_n$ is group-like, let $\ell := \log (X)$.  We
have $X = 1 + \ell + \ell^2/2! + \cdots $; here $\ell^k\in S^k(\widehat{\mathfrak f}_n)$,
so $\widehat p_n(X) = \ell$. Hence $\widehat{{\rm cbh}}_n(X) = \ell =
\log (X)$. 
\end{proof}

\section{Corollaries}

The KZ associator is def\/ined as follows.  Let $A_0$, $A_1$ be
noncommutative variables. Let $F_2$ be the free associative algebra
generated by $A_0$ and $A_1$, let ${\mathfrak f}_2 \subset F_2$ be its (free) Lie
subalgebra generated by $A_0$ and $A_1$. Let $\widehat F_2$ and $\widehat{\mathfrak f}_2$
be the degree completions of $F_2$ and ${\mathfrak f}_2$ ($A_0$ and $A_1$ have
degree $1$).

The KZ associator $\Phi$ is def\/ined \cite{Dr} as the renormalized
holonomy from $0$ to $1$ of the dif\/ferential equation
\begin{equation} \label{DE}
G'(z) = \left({{A_0}\over z} + {{A_1}\over{z-1}}\right) G(z) ,
\end{equation} 
i.e.,  $\Phi = G_1 G_0^{-1}$, where $G_0,G_1\in \widehat F_2 \, \widehat\otimes \,
{\mathcal O}_{]0,1[}$ are the solutions of (\ref{DE}) with $G_0(z) \sim z^{A_0}$
as $z\to 0^+$ and $G_1(z) \sim (1-z)^{A_1}$ as $z\to 1^-$; here
${\mathcal O}_{]0,1[}$ is the ring of analytic functions on $]0,1[$, and $\widehat
F_2 \,\widehat\otimes\, V$ is the completion of $F_2 \otimes V$ w.r.t. the
topology def\/ined by the $F_2^{\geq n} \otimes V$ (here $F_2^{\geq n}$
is the part of $F_2$ of degree $\geq n$).

We recall Le and Murakami's formula for $\Phi$ \cite{LM}.  We say
that a sequence $(a_1,\ldots,a_n)\in \{0,1\}^n$ is admissible if $a_1
=1$ and $a_n =0$. If $(a_1,\ldots,a_n)$ is admissible, we set 
\[
\omega_{a_1,\ldots,a_n} = \int_0^1 \omega_{a_1} \circ \cdots 
\circ \omega_{a_n} ,  
\]
where $\omega_0(t) = dt/t$, $\omega_1(t) = dt/(t-1)$ and $\int_a^b 
\alpha_1 \circ \cdots \circ \alpha_n = \int_{a\leq t_1 \leq \cdots 
\leq t_n \leq
b} \alpha_1(t_1) \wedge \cdots \wedge \alpha_n(t_n)$. Up to sign, the
$\omega_{a_1,\ldots,a_n}$ are MZV's (multiple zeta values).

If $(i_1,\ldots,i_n)$ is an arbitrary sequence in $\{0,1\}^n$, and
$(a_1,\ldots,a_n)$ is an admissible sequence, def\/ine integers
$C_{i_1,\ldots,i_n}^{a_1,\ldots,a_n}$ by the relation 
\begin{gather*}
\sum_{(i_1,\ldots,i_n) \in \{0,1\}^n} C_{i_1,\ldots,i_n}^{a_1,\ldots,a_n}
A_{i_n}\cdots A_{i_1}\\
\qquad
= \sum_{\mbox{\scriptsize $\begin{array}{c} S\subset \{\alpha | a_\alpha=0\},\\
T\subset \{\beta | a_\beta =1\} \end{array}$}}
(-1)^{{\rm card}(S) + {\rm card}(T)} A_1^{{\rm card}(T)} 
A(a_1,\ldots,a_n)^{S,T} A_0^{{\rm card}(S)},  
\end{gather*} 
where for any $S \subset \{\alpha | a_\alpha=0\}$, $T\subset \{\beta |
a_\beta = 1\}$, $A(a_1,\ldots,a_n)^{S,T} := \prod\limits_{\alpha\in [1,n]
\setminus(S\cup T)} A_{a_\alpha}$ (the product is taken in decreasing order of
the $\alpha$'s).

\begin{theorem}[\cite{LM}] \label{thm:1}
\[
\Phi = 1 + \sum_{n\geq 1} 
\sum_{\mbox{\scriptsize $\begin{array}{c} (a_1,\dots ,a_n){\rm \ admissible}\\
(i_1,\dots ,i_n)\in\{0,1\}^n\end{array}$}}
\omega_{a_1,\ldots,a_n}
C^{a_1,\ldots,a_n}_{i_1,\ldots,i_n} 
A_{i_n} \cdots A_{i_1}. 
\]
\end{theorem}

Since $\Phi\in\widehat F_2$ is a group-like element, Proposition~\ref{main} 
implies that $\log (\Phi) = \widehat{{\rm cbh}}_2(\Phi)$, 
therefore: 

\begin{corollary} \label{cor:phi}
\[
\log (\Phi) = \sum_{n\geq 1} 
\sum_{\mbox{\scriptsize $\begin{array}{c} (a_1,\dots ,a_n){\rm \ admissible}\\
(i_1,\dots ,i_n)\in\{0,1\}^n\end{array}$}}
\omega_{a_1,\ldots,a_n}
C^{a_1,\ldots,a_n}_{i_1,\ldots,i_n} 
{\rm CBH}_n(A_{i_n}, \ldots, A_{i_1}). 
\]
\end{corollary} 

Using the explicit formula of \cite{E}, one computes similarly 
the logarithm of the analogue $\Psi$ of the KZ associator of the equation
$G'(z) = \big(A/z + \sum\limits_{\zeta|\zeta^n = 1} b_{\zeta}/(z-\zeta)\big) G(z)$. 

Proposition \ref{main} also implies:  

\begin{lemma} \label{lemma}
Let ${\mathfrak g}$ be a nilpotent Lie algebra, $G$ be the associated Lie group, 
let $a<b\in {\mathbb R}$. 
Fix $h(z) \in C^0([a,b],{\mathfrak g})$ and let $H$ be the holonomy from
$a$ to $b$ of the differential equation $H'(z) = h(z) H(z)$, where $H(z)
\in C^1([a,b],G)$. Then  
\[
\log (H) = \sum_{n\geq 1} \int_{a\leq z_1\leq \cdots\leq z_n\leq b}
{\rm CBH}_n(h(z_n),\ldots,h(z_1)) dz_1 \cdots dz_n. 
\]
\end{lemma}

\subsection*{Acknowledgements} 

We f\/irst established the formula for 
$\log (\Phi)$ in Corollary~\ref{cor:phi} by analytic computations
(using a~direct proof of Lemma~\ref{lemma}). 
It was the referee who remarked its formal similarity with the 
formula of Le and Murakami (Theorem \ref{thm:1}); this remark can be 
expressed as the equality $\log (\Phi) = \widehat{{\rm cbh}}_2(\Phi)$. 
This led us to try and understand whether this formula followed
from the group-likeness of $\Phi$, which is indeed the case (Proposition~\ref{main}).
 C.~Reutenauer then pointed out that 
a part of our argument is a result in his book. 

\LastPageEnding

\end{document}